\documentclass[a4paper,12pt]{article}
\usepackage[english]{babel}
\usepackage[cp1250]{inputenc}
\usepackage{amssymb,amsmath}
\usepackage{enumerate}
\usepackage{graphicx}

\newcommand{\kart}{\, {\scriptstyle \square} \,}
\newcommand{\poln}{\, {\scriptstyle \boxtimes} \,}

\newtheorem{izrek}{Theorem}[section]
\newtheorem{posl}{Corollary}[section]
\newtheorem{lema}{Lemma}[section]
\newtheorem{primer}{Example}[section]

\newenvironment{dokaz}{\vskip6pt \noindent \textbf{Proof:}}{\hfill$\square$\\}

\title{Distance-residual graphs}
\author{Primož Lukšič\\ \small IMFM, University of Ljubljana, Ljubljana\\[-5pt] \small Slovenia\\
\small \texttt{primoz.luksic@fmf.uni-lj.si}\\[12pt]
 Tomaž Pisanski\\ \small IMFM, University of Ljubljana, Ljubljana\\[-5pt]
\small and University of Primorska, Koper\\[-5pt] \small Slovenia\\
\small \texttt{tomaz.pisanski@fmf.uni-lj.si}}

\date{September 15, 2006}

\begin{document}

\selectlanguage{english}

\maketitle
\begin{abstract}
If we are given a connected finite graph $G$ and a subset of its
vertices $V_{0}$, we define a distance-residual graph as a graph
induced on the set of vertices that have the maximal distance from
$V_{0}$. Some properties and examples of distance-residual graphs of
vertex-transitive, edge-transitive, bipartite and semisymmetric
graphs are shown. The relations between the distance-residual
graphs of product graphs and their factors are shown.\\[3pt]
\end{abstract}
\section{Introduction}

\noindent Let $G$ be a connected finite graph and let $V_{0} \subset
V(G)$ be a nonempty subset of vertices of $G$. We may form a
distance partition $P(G,V_{0})$ of $G$ with respect to $V_{0}$,
where $P(G,V_{0}) = \{V_0,V_1, \dots V_r\}$ and \begin{align*}
V_0 \cup V_1 \cup \dots \cup V_r & = V(G),\\
V_i \cap V_j & = \emptyset, \text{ for } i \neq j, \\
V_i & \neq \emptyset, \text{ for all } i.
\end{align*}

The classes $V_i$ are defined recursively as
\[V_{i+1} := \{v \in
V(G) \setminus (V_0 \cup V_1 \cup \dots \cup V_i)\: |\: \exists u
\in V_i \backepsilon: [u,v]\in E(G)\}.\]

The set $V_{i}$ contains all vertices of $G$ that have a minimal
distance of $i$ to the vertices of $V_{0}$ where the distance
$d(u,v)$ between two vertices $u$ and $v$ is defined as the shortest
path between them. Therefore $d(v,v_{i})\geq i$ for $v\in V_{0}$ and
$v_{i}\in V_{i}$, and there exists a vertex $v_{0}\in V_{0}$ for
which $d(v_{0},v_{i})= i$.

We are interested in induced subgraphs $\langle V_i \rangle$ defined
by the distance classes, particularly in the subgraph with the
lowest index $R_{G}:=\langle V_0 \rangle$, which we call the
\emph{root}, and the subgraph with highest index
$Res(G,R_{G}):=\langle V_r \rangle$, which we call the {\em
distance-residual graph} or the \emph{distance residual}. When the
root consists of a single vertex, i.e.\ $R_{G}\cong K_{1}$, the
residual is called a {\em vertex residual}. When $R_{G} \cong
K_{2}$, the residual is called an {\em edge residual}. With respect
to the definition of the distance residuals, all of the graphs in
the paper will be nontrivial, simple, finite, and in most cases
connected. Also some standard labels for some known graphs will be
used: $K_{n}$ for complete graphs, $K_{m,n}$ for complete bipartite
graphs, $C_{n}$ for cycles, $P_{n}$ for paths, and $m K_{n}$ for a
disjoint union of $m$ complete graphs on $n$ vertices.

The motivation for the definition of distance-residual graphs was in
extending the definition of \emph{distance sequence} which is an
ordered list where the $i$-th element equals the number of vertices
at distance $i$ from the selected root. The distance sequence
therefore presents only the number of vertices at a distance $i$
from the root but we are also interested in the induced subgraphs on
those vertices, especially on the set farthest away from the root.

In the next section we present some properties of the
distance-residual graphs with the focus on the vertex- and
edge-transitive graphs \cite{babai1,babai2}, bipartite graphs, and
semisymmetric graphs \cite{semisym1,mar}. In section 3 we show how
the distance residuals of product graphs for some well-known
products depend on the distance residuals of their factors. We
conclude with some open questions regarding distance residuals and
other distance related problems.

\section{Properties}

We can only define a distance-residual graph of a connected graph
but any graph, connected or not, can be a distance-residual graph.

\begin{izrek}
Let $H$ be an arbitrary graph and $n\in \mathbb{N}$. Then there
exists a connected graph $G$ with the root $R_{G}$ of order $n$ such
that $H$ is isomorphic to $Res(G,R_{G})$.
\end{izrek}
\begin{dokaz} Let us choose an arbitrary graph $R_{G}$
of order $n$ with the property $V(R_{G}) \cap V(H) = \emptyset$.
Graph $G$ is constructed as follows: The vertex set $V(G)$ consists
of $V(H) \cup V(R_{G})$. There are two types of edges in $G$. All
the original edges of $H$ and $R_{G}$ remain edges in $G$ and for
each vertex $v \in V(H)$ and each vertex $r \in V(R_{G})$ there is
an edge between them. Clearly, the distance partition is given by
$V_0 = V(R_{G})$ and $V_1 = V(H)$, therefore $Res(G,R_{G}) = H.$
\end{dokaz}

Most of the interesting cases occur for vertex-transitive and
edge-transitive graphs.

\subsection{Vertex- and edge-transitive graphs}

\begin{lema}\label{v-tranz}
Let $G$ be a connected vertex-transitive graph. Then all of its
vertex residuals are isomorphic, i.e.\ they do not depend on the
choice of the vertex for the root. The converse is not true even if
the graph is regular.
\end{lema}

\begin{dokaz}
The first part is obvious because if the vertex residuals would not
be isomorphic, the graph would not have a transitive automorphism
group. The graph in Figure \ref{fig:graph1} proves that the converse
is not true. It is built from two copies of $K_{4}$ by joining some
of their their edges (see \cite{pet} for a definition of joining)
and is therefore 3-regular. All of its vertex residuals are
isomorphic (to $K_{1}$) but it is not vertex-transitive because the
automorphism, which would map a vertex from $K_{4}$ to a vertex
created by joining edges, does not exist.
\begin{figure}[htb!]
\centering
\includegraphics[width=0.4\textwidth]{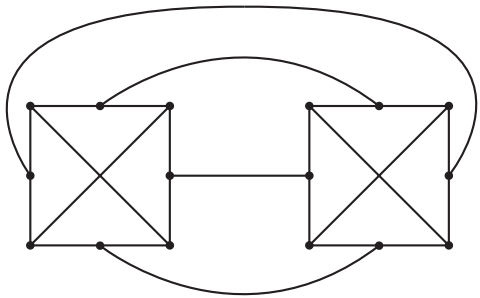}
\caption{Non-vertex-transitive cubic graph with isomorphic vertex
residuals.} \label{fig:graph1}
\end{figure}

\end{dokaz}

Vertex residual is also connected to the \emph{growth} of the graph
(see \cite {pis}) because its order is the leading coefficient of
the growth polynomial of the graph at its root. This is equivalent
to taking the last element of the distance sequence of the root.

We know that there exist graphs that are \emph{growth-regular},
i.e.\ their growth function is independent of the root vertex, but
not vertex-transitive. From the above proof we can also see that
such graphs can have isomorphic vertex residuals and still not be
vertex-transitive. The growth function of the graph in Figure
\ref{fig:graph1} is namely $1+3x+6x^2 +5x^3 + x^4$ independent of
the vertex we take for the root.

\begin{lema}
Let $G$ be a connected edge-transitive graph. Then all of its edge
residuals are isomorphic, i.e.\ they do not depend on the choice of
the two adjacent vertices for the root. The converse is not true.
\end{lema}

\begin{dokaz}
The first part is similar as in the prior lemma. As a counterexample
of the converse we present the graph in Figure \ref{fig:graph2} that
has all of its edge residuals isomorphic (to $K_{1}$) but it is
clearly not edge-transitive.

\begin{figure}[htb!]
\centering
\includegraphics[width=0.4\textwidth]{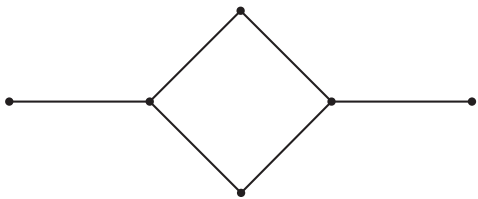}
\caption{Non-edge-transitive graph with isomorphic edge residuals.}
\label{fig:graph2}
\end{figure}
\end{dokaz}

\begin{lema}
Let $G$ be a connected edge-transitive graph and $L(G)$ its line
graph. Then all of the distance residuals $Res(L(G),K_{1})$ are
isomorphic.
\end{lema}
\begin{dokaz}
The line graph of a connected edge-transitive graph is a connected
vertex-transitive graph. The rest follows from Lemma \ref{v-tranz}.
\end{dokaz}

We present the distance residuals of some well-known graphs.

\begin{primer}
\begin{align*}
Res(C_{2k},K_{1}) &\cong K_{1}\quad & Res(C_{2k},K_{2})  &\cong  K_{2}\\
Res(C_{2k+1},K_{1}) & \cong K_{2}\quad & Res(C_{2k+1},K_{2}) &\cong  K_{1}\\
Res(K_{n},K_{1}) &\cong K_{n-1} \quad  & Res(K_{n},K_{k}) &\cong K_{n-k}\\
\end{align*}
\end{primer}

\begin{primer}
The Petersen graph $P(5,2)$ is vertex-transitive, edge-transitive
and also distance-transitive, therefore it has isomorphic distance
residuals for some roots:
\begin{align*}
Res(P(5,2),K_{1}) &\cong C_{6},\\
Res(P(5,2),K_{2}) &\cong  2 K_{2},\\
Res(P(5,2),P_{3}) &\cong  2 K_{1},\\
Res(P(5,2),P_{4}) &\cong  K_{1},\\
Res(P(5,2),C_{5}) &\cong  C_{5}.
\end{align*}
\end{primer}

\subsection{Bipartite graphs}

For bipartite graphs a special relation between vertex and edge
residuals holds. But before we present it, we shall define the
distance between vertices and subgraphs. Let $d(R,v)$ be the
distance from $R\subset G$ to vertex $v\in V(G)$, i.e.\
$d(R,v):=\min_{r\in V(R)} d(r,v)$, and let $d(R,H):=\min_{v\in
V(H)}d (R,v)$ be the distance between two subgraphs of graph $G$. If
graph $G$ is connected, the distances are well defined.

We will often need the distance between the root and the distance
residual of some connected graph $G$, so we will denote
$d(R_{G},Res(G,R_{G}))$ as $d_{R_{G}}$. The distance namely depends
on the root of the graph.

In the next theorem we use the above notation to shorten the writing
for the distance $d(\langle \{u\}\rangle,Res(G,\langle
\{u\}\rangle))$ to $d_{u}$ and similarly to shorten $d(\langle
\{u,v\}\rangle,Res(G,\langle \{u,v\}\rangle))$ to $d_{u,v}$.

\begin{izrek} \label{bipart}
Let $G \not\cong K_{2}$ be a connected bipartite graph with
partitions $P_{1}$ and $P_{2}$. Let $u\in P_{1}$ and $v\in P_{2}$ be
two neighbors from $G$. Then the edge residual equals
\[Res(G,\langle \{u,v\}\rangle) =\left\{
  \begin{array}{ll}
    \left\langle \, V(Res(G,\langle \{u\}\rangle)) \cup  V(Res(G,\langle \{v\}\rangle))
  \, \right\rangle, & \hbox{if } d_{u}=d_{v};\\
    Res(G,\langle \{u\}\rangle), & \hbox{if } d_{u}=d_{v}+1;\\
    Res(G,\langle \{v\}\rangle), & \hbox{if } d_{v}=d_{u}+1.
  \end{array}
\right.
\]
There are no other possibilities.
\end{izrek}

\begin{dokaz}
First we prove that the three cases above are the only ones
possible. Let us assume that $d_{u}>d_{v}+1$ and take a vertex $r\in
V(Res(G,\langle \{u\}\rangle))$. From $d_{v}\geq d(v,r)$ we get
$d_{u}> d(v,r)+1$. But that means there is a shorter path between
$u$ and $r$ via $v$ ($u$ and $v$ are neighbors) which is a
contradiction. The proof is similar if we assume $d_{v}>d_{u}+1$.
From this we can also deduce that $d(u,s)=d(v,s)\pm 1$ for any
vertex $s\in V(G)$; because the graph is bipartite, the distances
cannot be equal.

If $d_{u}$ and $d_{v}$ are equal, the vertex residuals lie in
different partitions. If we take an arbitrary vertex $r\in
V(Res(G,\langle \{u\}\rangle))$, then $d(v,r)=d_{u}-1$ (if
$d(v,r)=d_{u}+1$, then $d_{v}>d_{u}$). The same is true if we change
the role of $u$ and $v$. So $d_{u,v}\geq d_{u}-1=d_{v}-1$.

To prove that $d_{u,v}$ is equal to $d_{u}-1$ and that the
edge-residual graph is induced only on the vertices of both vertex
residuals we take a vertex $s$ which is not in any of the vertex
residuals. If $d_{u}> d(u,s)+1$, then $d_{u}-1>d(u,s)$. Let $d_{u}=
d(u,s)+1$. Then $d(v,s)=d(u,s)-1$ because if
$d(v,s)=d(u,s)+1=d_{u}=d_{v}$, $s$ would be in the vertex residual
of $v$. It follows that $d(\langle \{u,v\}\rangle
,s)=d(u,s)-1=d_{u}-2<d_{u}-1$. The proof is the same if we
interchange $u$ and $v$.

If $d_{u}=d_{v}+1$ and we take a vertex $r\in V(Res(G,\langle
\{u\}\rangle))$, then $d_{u}=d(u,r)=d(v,r)+1$, therefore
$d(v,r)=d_{v}$ and $r\in V(Res(G,\langle \{v\}\rangle))$. It follows
that $d_{u,v}=d_{v}$. If we take a vertex $s$ from the vertex
residual of $v$ but not of $u$, then $d_{v}=d(v,s)=d(u,s)+1$.
Therefore $d(\{u,v\},s)=d(u,s)=d(v,s)-1<d_{v}$ and $s$ is not in the
edge residual. The proof in the third case is essentially the same;
we just interchange $u$ and $v$.

\end{dokaz}

\begin{primer}
Let $K_{m,n}$ ($m,n\neq 1$) have partitions sets $|P_{1}|=m$ and
$|P_{2}|=n$. Let $u\in P_{1}$ and $v\in P_{2}$. Then
$Res(K_{m,n},\langle \{u\}\rangle)) \cong (m-1)K_{1}$ and
$Res(K_{m,n},\langle \{v\}\rangle))\cong (n-1)K_{1}$. Because
$d_{u}=d_{v}=2$, we get from Theorem \ref{bipart}
$Res(K_{m,n},K_{2}) \cong (m+n-2)K_{1}$.
\end{primer}

An interesting class of bipartite graphs are \emph{semisymmetric
graphs}, i.e.\ regular graphs which are edge-transitive but not
vertex-transitive. Semisymmetric graphs have an automorphism group
that acts transitively on each of the bipartition sets. That means
the vertex residuals of roots from the same set are isomorphic.
Furthermore, the distance sequences in each of the sets are also the
same.

In fact, one motivation for the definition of distance-residual
graphs was the fact that the \emph{Gray graph} \cite{gray2, gray4},
the smallest semisymmetric cubic graph \cite{gray1}, is the edge
residual of the generalized quadrangle $W(3)$ \cite{quad2, quad} as
mentioned in \cite{gray3}. Gray graph is of order 54 and has the
distance sequences $(1,3,6,12,12,12,8)$ and $(1,3,6,12,16,12,4)$
with the vertex residuals $8K_{1}$ and $4K_{1}$ depending on the
partition from which the root vertex was taken. The edge residual is
induced on the 12 vertices of both vertex residuals and is
isomorphic to $4P_{3}$.

We mention in passing that there is an error in \cite{gray2} in the
Figure 7 which represents the construction 2.4 of the Gray graph.
The vertex labels $a_{i}$ and $b_{i}$ (for $1\leq i \leq 4$) in the
auxiliary graph $G_{2345}$ must me interchanged.

The smallest semisymmetric graph is the \emph{Folkman graph} of
order 20 and valence 4 which has distance sequences $(1, 4, 9, 6)$
and $(1, 4, 6, 6, 3)$, and vertex residuals $6K_{1}$ and $3K_{1}$.
By Theorem \ref{bipart} the edge residual graph equals $3K_{1}$.

For another example we look at the so called \emph{Ljubljana graph}
\cite{ljubljana} of order 112 with the distance sequences $(1, 3, 6,
12, 24, 34, 24, 7, 1)$ and $(1, 3, 6, 12, 24, 34, 25, 7)$. The
vertex residuals are isomorphic to $K_{1}$ and $7 K_{1}$ with the
edge residual equal to $K_{1}$. The Ljubljana graph was originally
discovered by R.~Foster (unpublished) and later studied in a series
of papers by I.J.~Dejter and his co-authors \cite{ljubljana3,
ljubljana1, ljubljana2, ljubljana4}. Only in \cite{ljubljana}, where
the Ljubljana graph was rediscovered for the third time, it was
determined that it is the unique third smallest cubic semisymmetric
graph, and hence isomorphic to the graph of Foster, Dejter, et al.

\section{Distance residuals of product graphs}

Product graphs have various interesting properties that make them
subject of intensive studies. Problems, that are intractable for
general graphs, sometimes admit elegant solutions for special
classes of graphs, such as product graphs. This fact frequently drew
our attention in the past \cite{product6, product1, product2,
product4, product7, product5, product3}.

In this section we show how distance-residual graphs of product
graphs depend on the factors of those products and on their
respective distance-residual graphs. All of the well-known graph
products are covered: Cartesian, strong, direct, and lexicographic
product (see \cite{product} for more about graph products and their
properties). This simplifies the discovery of distance-residual
graphs in some well-known graphs as well as proving some interesting
properties of vertex-transitive graphs.

\subsection{Cartesian product}

We start with the Cartesian product (denoted as $G \kart H$), which
is the most fundamental and the most studied of all. Its vertex set
is, like the sets of all the other here mentioned products, defined
on the Cartesian product $V(G)\times V(H)$ of the vertex sets of the
factors. Its edge set is the set of all pairs $[(u,x),(v,y)]$ where
either $u=v$ and $[x,y]\in E(H)$ or $x=y$ and $[u,v]\in E(G)$.

The distance-residual graph of the Cartesian product of two graphs
can easily be expressed with the distance-residual graphs of
the respective factors.\\

\begin{izrek}
Let $G$ and $H$ be two connected graphs with $R_{G}\subset G$ and
$R_{H}\subset H$ as their roots. We can state the following
connection between distance-residual graphs:
\[ Res(G \kart H, R_{G}\kart R_{H})\cong Res(G,R_{G})\kart Res(H,R_{H}).
\]
\end{izrek}

\begin{dokaz}
Because both factors are connected, so is the product graph and
therefore the distance-residual graph $Res(G \kart H, R_{G}\kart
R_{H})$ is well defined. The distances in the product are sums of
distances in both factors \cite{product}, i.e.\ the distance between
two vertices $(g_{1},h_{1})$ and $(g_{2},h_{2})$ in $G \kart H$ is
equal to $d_{G}(g_{1},g_{2})+d_{H}(h_{1},h_{2})$. Therefore
\[
d_{G \kart H}(R_{G}\kart R_{H},(g,h))=d_{G}(R_{G},g)+d_{H}(R_{H},h),
\]
because \begin{align*}\min_{(r_{G},r_{H})\in V(R_{G}\kart
R_{H})}d_{G \kart H}((r_{G},r_{H}),(g,h)) &=\min_{(r_{G},r_{H})\in
V(R_{G}\kart R_{H}(} d_{G}(r_{G},g)+d_{H}(r_{H},h)\\ &=
\min_{r_{G}\in V(R_{G})} d_{G}(r_{G},g)+  \min_{r_{H}\in
V(R_{H})}d_{H}(r_{H},h)\\ & = d(R_{G},g)+d(R_{H},h).\end{align*}

So the residual graph contains the vertex $(g,h)$ if and only if
$g\in V(Res(G,R_{G}))$ and $h\in V(Res(H,R_{H}))$. And finally,
because the distance-residual graph is induced, we must take the
Cartesian product of distance-residual graphs of both factors.
\end{dokaz}

From the associativity of the Cartesian product it also follows:
\begin{posl}
\[ Res(G_{1} \kart \cdots \kart G_{n}, R_{G_{1}}\kart \cdots \kart R_{G_{n}}) \cong
Res(G_{1},R_{G_{1}})\kart \cdots \kart Res(G_{n},R_{G_{n}}).
\]
\end{posl}

\begin{primer}
The $n$-dimensional hypercube $Q_{n}$ is the Cartesian product of
$n$ copies of $K_{2}$, so its vertex residual is $K_{1}$. Also, if
some graph in the product has a trivial vertex residual, we do not
have to take it into consideration.
\end{primer}

\subsection{Strong product}

Strong product, which is denoted by $G \poln H$, is similar to the
Cartesian product because the edge set of the product includes the
same edges as the Cartesian product with the addition of those edges
$[(u,x),(v,y)]$ where $[u,v]\in E(G)$ and $[x,y]\in E(H)$.

There are three possible options for the distance-residual graph of
the strong product, depending on the distances from roots to the
points in distance-residual graphs in both factors. So we reuse the
notation $d_{R_{G}}$ as the distance $d(R_{G},Res(G,R_{G}))$ from
the root to the distance-residual graph of some graph $G$.

\begin{izrek}
Let $G$ and $H$ be two connected graphs with $R_{G}\subset G$ and
$R_{H}\subset H$ as their roots. The distance-residual graph of the
strong product of these graphs is
\[Res(G \poln H, R_{G}\poln R_{H})\cong
\left\{
  \begin{array}{ll}
    H \poln Res(G,R_{G}), & \hbox{if } d_{R_{G}}> d_{R_{H}};\\
    G \poln Res(H,R_{H}), & \hbox{if } d_{R_{G}}< d_{R_{H}};\\
    \langle \, V(Res(G,R_{G}))\times V(H) \:\cup &\\
    \;\:V(G) \times V(Res(H,R_{H})) \, \rangle, & \hbox{else.}
  \end{array}
\right.
\]
\end{izrek}

\begin{dokaz}
We follow the proof presented in the case of the Cartesian product.
Here also the graph of the strong product is connected so the
distance-residual graph is well defined. The distance in the product
is equal to the maximal distance in both factors \cite{product},
therefore
\[
d_{G \poln H}(R_{G}\poln R_{H},(g,h))=\max \{
d_{G}(R_{G},g),d_{H}(R_{H},h)\}.
\]

If $d_{R_{G}}>d_{R_{H}}$, it follows $d_{G \poln H}(R_{G}\poln
R_{H},(g,h))= d(R_{G},g)$. So the distance-residual graph contains
the vertex $(g,h)$ if and only if $g\in V(Res(G,R_{G}))$. Vertex $h$
is therefore an arbitrary vertex from graph $H$ and because the
distance-residual graph is induced, we get the mentioned result. The
same argument follows when $d_{R_{G}}< d_{R_{H}}$.

In the case where both distances are equal we get a
distance-residual graph which comprises of all the vertices of
$Res(G,R_{G})\times H$ and $G \times Res(H,R_{H})$ (with the
vertices of $Res(G,R_{G}) \times Res(H,R_{H})$ represented twice).
Once again, because the distance-residual graph is induced, we get
the mentioned result.
\end{dokaz}

Let $\mathcal{G}_{\max}$ be the set of all connected graphs
$G_{1},\ldots , G_{n}$ with their respective roots $R_{G_{1}},\ldots
, R_{G_{n}}$ for which the distance $d_{R_{G_{i}}}$ is maximal. To
put it in another way, $\mathcal{G}_{\max}:=\left\{ G_{j} \, | \,
d_{R_{G_{j}}}= \max_{1\leq i \leq n} \{d_{R_{G_{i}}} \} \right\}$.
Let us rearrange the graphs so that the members of
$\mathcal{G}_{\max}$ get the indexes from one to $k\leq n$. In a
similar way as with the Cartesian product we can come to the
following conclusion:
\begin{posl}
If $\mathcal{G}_{\max}$ has $k$ $(1\leq k \leq n)$ elements, then
the distance-residual graph $Res(G_{1}\poln \cdots \poln
G_{n},R_{G_{1}}\poln \cdots \poln R_{G_{n}})$ is isomorphic to
\begin{equation}  \left\langle \, \bigcup_{1\leq i \leq k} V(G_{1}) \times
\ldots \times V(G_{i-1}) \times V(Res(G_{i},R_{G_{i}})) \times
V(G_{i+1}) \times \ldots \times V(G_{n}) \,
\right\rangle.\label{strong-res} \end{equation} If
$\mathcal{G}_{\max}$ has only one member ($G_{1}$), the graph in
(\ref{strong-res}) can be written as
\[Res(G_{1},R_{G_{1}})\poln G_{2} \poln \cdots \poln G_{n}.\]
\end{posl}

\begin{primer}
The complete graph on $m n$ vertices is the strong product of
$K_{m}$ with $K_{n}$. Therefore, $Res(K_{mn},K_{s}\poln K_{t})\cong
K_{mn-st}$ for $s<m$ and $t<n$.
\end{primer}

\subsection{Lexicographic product}

The major difference between the lexicographic product (denoted here
as $G \circ H$) and the products mentioned before is that this
product is not commutative. The lack of this symmetry can be seen
from the definition of edges because $[(u,x),(v,y)]\in E(G \circ H)$
if either $[u,v]\in E(G)$ or $u=v$ and $[x,y]\in E(H)$.

The product is sometimes written as $G[H]$ and can be represented by
replacing each vertex $v\in V(G)$ with the copy of graph $H$,
denoted $H_{v}$, and then connecting all off the vertices of $H_{u}$
with all off the vertices of $H_{v}$ if the vertices $u$ and $v$ are
neighbors in $G$. The distance between vertices in the product
therefore depends on whether they lie in the same copy of $H$. If
they do, they are at a distance two (or one if they are neighbors),
if they do not, their distance is determined by the distance of
their first factors, i.e.\ the distance between the copies of $H$ in
which they lie in. Therefore $d_{R_{G\circ H}}$ can be equal to 1, 2
or $d_{R_{G}}$.

From the above description we can see that if $G$ is connected (and
nontrivial) then so is $G\circ H$ for arbitrary $H$. The statement
is true also in the other direction (see \cite{product}).

The same as with the strong product, the distance-residual graph
depends on the distances between roots and distance-residual graphs
of its factors. Therefore, we use the notation established there.

\begin{izrek} \label{lex}
Let $G$ and $H$ be two graphs with $R_{G}\subset G$ and
$R_{H}\subset H$ as their roots and let $G$ be connected. The
distance-residual graph of the lexicographic product of these
graphs, $Res(G \circ H, R_{G}\circ R_{H})$, is isomorphic to
\begin{enumerate}
  \item $\displaystyle \bigcup_{\substack{r_{G}\in V(R_{G}) \\ r_{G}\; isolated}}\!\!\!\!\! \left(H \setminus \langle V_{0} \cup V_{1}\rangle_{H}\right)$, if
$d_{R_{G}}=1$, $d_{R_{H}}\neq 1$ and there exist isolated vertices
in $R_{G}$,
  \item $(G\circ H) \setminus (R_{G}\circ R_{H})$, if $d_{R_{G}}=1$ and ($d_{R_{H}}=1$ or $R_{G}$ has no
isolated vertices),
  \item $Res(G,R_{G})\circ H \:\cup \!\! \displaystyle
\bigcup_{\substack{r_{G}\in V(R_{G}) \\ r_{G}\; isolated}}\!\!\!\!\!
\left(H \setminus \langle V_{0} \cup V_{1}\rangle_{H}\right)$, if
$d_{R_{G}}=2$,
  \item $Res(G,R_{G})\circ H$, if $d_{R_{G}}\geq 3.$
\end{enumerate}
\end{izrek}

\begin{dokaz}
\begin{enumerate}[1.]
  \item If $d_{R_{G}}=1$, then all of the vertices in $H_{v}$, where $v\not\in
V(R_{G})$, lie in the class $V_{1}$ of the distance partition of
$G\circ H$. If two vertices $u$ and $v$ from $R_{G}$ are adjacent,
then the vertices of $V(H_{v})\setminus V(R_{H_{v}})$ and
$V(H_{u})\setminus V(R_{H_{u}})$ also lie in $V_{1}$. But if $v$ is
an isolated vertex of $R_{H}$, i.e.\ its neighbors are from
$V(G)\setminus V(R_{G})$, then the vertices of $H_{v} \setminus
\langle V_{0} \cup V_{1}\rangle_{H_{v}}$ (which exist because
$d_{H}\neq 1$) lie in $V_{2}$ and therefore induce the
distance-residual graph of the product.
  \item If $d_{R_{G}}=1$ and $R_{G}$ has no isolated vertices, then, following the above
proof, we see that all of the vertices, which are not from
$R_{G}\circ R_{H}$, lie in $V_{1}$. The same follows if
$d_{R_{H}}=1$.
  \item If $d_{R_{G}}=2$, then for all $v\in V(Res(G,R_{G}))$ the
vertices of $H_{v}$ lie in $V_{2}$ of the partition $P(G\circ
H,V_{0})$. But from the proof of the first case we see that $V_{2}$
also consists of vertices from $H_{v} \setminus \langle V_{0} \cup
V_{1}\rangle_{H}$ for isolated $v\in V(R_{G})$. If there are no
isolated vertices in $R_{G}$ or if $d_{R_{H}}=1$, then the latter
vertices do not exist.
  \item This follows from the third case because for all $v\in V(Res(G,R_{G}))$ the
vertices of $H_{v}$ lie in $V_{d_{R_{G}}}$ of the partition
$P(G\circ H,V_{0})$ which induces the distance-residual graph.
\end{enumerate}
\end{dokaz}

We can generalize the result in case 4 of the Theorem \ref{lex} to
more factors with the help of the associativity of the lexicographic
product and calculate it recursively.

\begin{posl}
Let $G_{1}, G_{2}, \ldots , G_{n}$ be graphs with their respective
roots $R_{G_{1}}, R_{G_{2}},  \ldots , R_{G_{n}}$ and let $G_{1}$ be
connected. Also, let the distance $d_{R_{G_{1}} \circ \cdots \circ
R_{G_{i}}}$ in $G_{1} \circ \cdots \circ  G_{i}$ be at least 3 for
all $i\in \{1,\ldots n-1\}$. It follows that
\[
Res(G_{1}\circ G_{2} \circ \cdots \circ G_{n}, R_{G_{1}} \circ
R_{G_{2}} \circ \cdots \circ R_{G_{n}}) \cong
Res(G_{1},R_{G_{1}})\circ G_{2} \circ \cdots \circ G_{n}.
\]
\end{posl}

With the use of the lexicographic product we can prove another
property of the distance-residual graphs in vertex-transitive
graphs.

\begin{izrek} \label{v-tranz:res}
Let $H$ be a vertex-transitive graph and $n\in \mathbb{N}$. Then
there exists a connected vertex-transitive graph $G$ with the root
$R_{G}$ of order $n$ such that $H$ is isomorphic to $Res(G,R_{G})$.
\end{izrek}

\begin{dokaz}
For the construction of graph $G$ we use the lexicographic product
of graphs $C_{k}$ and $H$. Because both of them are
vertex-transitive so is their product $G$ (see \cite{product}) and
because $C_{k}$ is connected so is the product. Let $R_{H}:=K_{1}$
and $R_{C_{k}}:=P_{n}$. Therefore $|R_{H}\circ R_{C_{k}}|=n$. If we
take $k=n+5$, then $d_{R_{C_{k}}}=3$. The $Res(C_{n+5} \circ H,
R_{C_{n+5}}\circ R_{H})$ is by case 4 of Theorem \ref{lex}
isomorphic to $Res(C_{n+5},P_{n})\circ H\cong K_{1}\circ H\cong H$.
\end{dokaz}

\begin{primer}
With the help of Theorem \ref{v-tranz:res} we can construct a
vertex-transitive graph that has all of its vertex residuals
isomorphic to the Petersen graph. But the latter is also the vertex
residual of the Clebsch graph shown on Figure \ref{fig:clebsch},
which is a symmetric (vertex- and edge-transitive), strongly regular
graph on 16 vertices of valence 5.
\begin{figure}[htb!]
\centering
\includegraphics[width=0.4\textwidth]{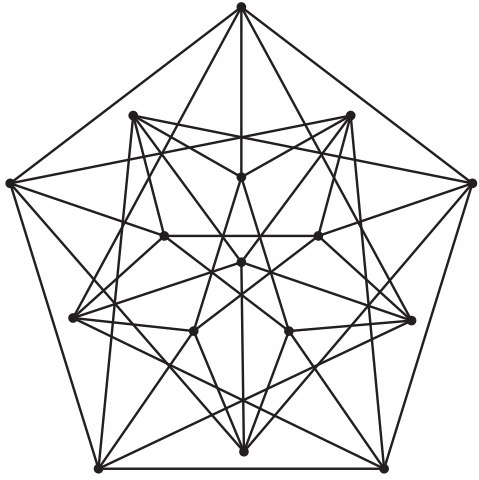}
\caption{The Clebsch graph which has all of its vertex residuals
isomorphic to the Petersen graph.} \label{fig:clebsch}
\end{figure}
\end{primer}


\subsection{Direct product}

Finally, we mention the direct product which is also known as tensor
product or categorical product \cite{product}. We denote it by $G
\times H$ (the sign $\otimes$ is also frequently used). Its edge set
is made up of edges $[(u,x),(v,y)]$ where $[u,v]\in E(G)$ and
$[x,y]\in E(H)$.

The product has some unusual properties. First of all, the
connectivity of both factors is not sufficient condition for the
connectivity of the product; at least one of them must not be
bipartite. Furthermore, the distance function between vertices in
the product is unlike with the other products, as shown by the
following lemma.

\begin{lema}[see \cite{tensor}] \label{tensor}
Let $x = (x_1, x_2,\ldots, x_{n})$ and $y = (y_1, y_2, \ldots ,
y_n)$ be two vertices of $G :=G_{1} \times G_{2} \times \ldots
\times G_{n}$. If there is no integer $m$ for which each $G_i$ has
an $x_{i}-y_{i}$ walk of length $m$, then $d_{G}(x, y) = \infty$.
Otherwise, $d_{G}(x, y) = \min \{m \in \mathbb{N}\,\,\, | \text{
each } G_{i} \text{ has an } x_{i}-y_{i} \text{ walk of length }
m\}$.
\end{lema}

We present a theorem that follows from Lemma \ref{tensor} from which
we can conclude that the distance-residual graph of the direct
product does not depend on the residuals of the factors.

\begin{izrek}
Let $G$ and $H$ be connected graphs and at least one of them
non-bipartite. Let $R_{G}$ and $R_{H}$ be their respective roots.
Then the distance-residual graph $Res(G\times H, R_{G}\times R_{H})$
is induced on all of the vertices $(u,v)\in V(G\times H)$ for which
the following holds:
\begin{align*}
d(R_{G}\times R_{H},(u,v)) =& \max_{(g,h)\in V(G\times H)}
\min_{\,\, r_{G}\in R_{G},r_{H}\in R_{H}} \{ m \in \mathbb{N}\,\,\,
| \,\, G \text{ has an } r_{G}-g \text{ walk} \\ &  \text{of length
} m \text{ and } H \text{ has an } r_{H}-h \text{ walk of length }
m\}.
\end{align*}

\end{izrek}

\begin{primer} \label{direkt}
Let $G:=P_{3}$ with $R_{G}$ consisting of one of the two noncentral
vertices of $P_{3}$. Let $H$ be a graph consisting of $K_{3}$ with
an appended edge on each of the vertices. For the root of $H$ we
take one of the vertices not from $K_{3}$. Then $Res(P_{3}\times
H,R_{G}\times R_{H})\cong K_{1}$ with the vertex having as the first
component the central vertex of $P_{3}$ and the second from the root
of $H$.
\begin{figure}[htb!]
\centering
\includegraphics[width=0.7\textwidth]{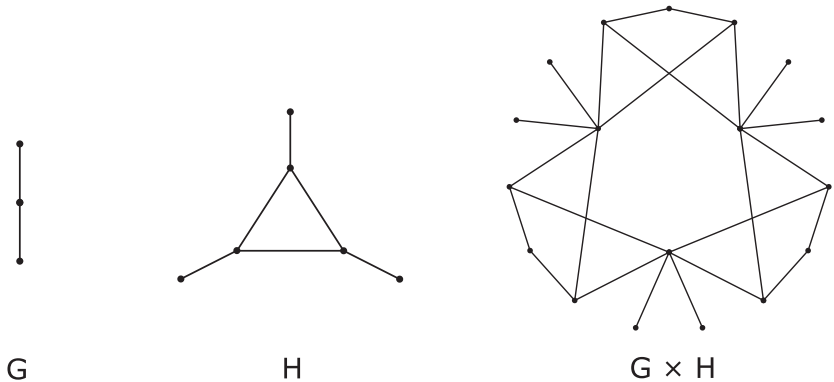}
\caption{Graphs $G$ and $H$ from Example \ref{direkt} and their
direct product}
\end{figure}
\end{primer}

\section{Concluding remarks}

In this paper we have defined a distance-residual graph and proven
some of its properties. It would be interesting to see whether our
methods could be used to answer some of the following questions.

Since every graph can be a distance-residual graph, it is a
challenge to find well-known graphs as distance residuals of some
other well-known graphs. We also ask what is the sufficient
condition for a growth-regular graph, i.e.\ a graph with the same
distance sequences, to be vertex-transitive. Graph bundles
\cite{bundle4, bundle1, bundle2, bundle3, bundle5} form an
interesting generalizations of product graphs. It would be of
interest to investigate their properties in connection to residual
graphs. And finally, regular edge-transitive graph that admits two
distinct distance sequences is necessarily semi-symmetric. In
principle, the converse need not be true. It would be interesting to
apply our methods for construction of families of semi-symmetric
graphs in which all vertices give rise to the same distance
sequence.

\end{document}